\documentclass[11pt]{article}

\usepackage{latexsym,epsfig} %% Add other packages as necessary
\usepackage{amsmath}
\usepackage{amsthm, amssymb}
\usepackage{mathtools}
\usepackage{graphicx}   % need for figures
\usepackage{url}
\usepackage[utf8]{inputenc} % Swedish
\usepackage[affil-it]{authblk}
\usepackage{mathrsfs} % for mathscr

\title{The Limit Shape of a Stochastic Bulgarian Solitaire}

\author[1]{Kimmo Eriksson}
\author[1]{Markus Jonsson}
\author[2]{Jonas Sjöstrand}
\affil[1]{Mälardalen University, School of Education, Culture and Communication, Box 883, SE-72123 Västerås, Sweden}
\affil[2]{Royal Institute of Technology, Department of Mathematics, SE-10044 Stockholm, Sweden}

\newtheorem{theorem}{Theorem}
\newtheorem{lemma}{Lemma}

\newtheorem{corollary}{Corollary}
\newtheorem{observation}{Observation}

\DeclareMathOperator{\ord}{\textup{\textsf{ord}}}

\newcommand{\Par}{\mathcal{P}}
\newcommand{\Bin}{\textup{Bin}}
\newcommand{\rescaled}[2]{\partial^{#1}{#2}}
\newcommand{\Rnn}{\mathbb{R}_{\geq 0}}

\begin{document}
\maketitle

\begin{abstract}
We consider a stochastic version of Bulgarian solitaire: A number of cards are distributed in piles; in every move a new pile is formed by cards from the old piles, and each card is picked independently with a fixed probability. This game corresponds to a multi-square birth-and-death process on Young diagrams of integer partitions. We prove that this process converges, in a strong sense, to an exponential limit shape as the number of cards tends to infinity. Furthermore, we bound the probability of deviation from the limit shape and relate this to the number of moves played in the solitaire.
\end{abstract}

\section{Introduction}
\label{sec:bulgsol}
The game of Bulgarian solitaire is played with a deck of $n$ identical cards divided arbitrarily into several piles. A move consists of picking one card from each pile and letting these cards form a new pile. This move is repeated over and over again. If the total number of cards in the deck is a triangular number, i.e., $n=1+2+\ldots+k$ for some $k$, this process has an interesting property: Regardless of the initial configuration, a finite number of moves will lead to the stable configuration where there is one pile of each size from 1 up to $k$, see \cite{Brandt, Vassilev}. This property motivated initial interest in Bulgarian solitaire in the early 1980's, and it  featured in a 1983 column by Martin Gardner in Scientific American. Later research have studied also other aspects of Bulgarian solitaire \cite{AkinDavis1985, Bentz1987, Etienne1991, Griggs1998, Igusa1985}.  For  information about the earlier history of the game (including its name), and a summary of subsequent research, see \cite{drensky2015bulgarian, Hopkins2011}.

If $n$ is not a triangular number, a stable configuration does not exist but after at most $O(n)$ moves the game will enter into a cycle. Moreover, all configurations of the cycle are ``almost triangular'' in the following sense: If $k=\max\{ j:1+2+\ldots+j\leq n \}$, then all the configurations in that cycle can be constructed from the triangular configuration $(k,k-1,\ldots,1)$ by adding at most one card to each pile, and possibly adding one more pile of size 1.  For exact formulations and more details, see \cite{AkinDavis1985,Bentz1987, Etienne1991, Griggs1998}. Obviously, in any configuration the set of pile sizes constitute an integer partition of $n$. Thus, as a sweeping statement we can say that the Bulgarian solitaire converges to a triangular shape of the corresponding Young diagram.

Popov \cite{Popov} considered a stochastic version of Bulgarian solitaire. In this version, a move in the game consists of forming a new pile by picking one card from each pile in a \emph{random selection} of the piles. Specifically, any given pile must release one card with a fixed probability $0<p<1$ and independently of the other piles. The resulting game is a discrete-time irreducible and aperiodic Markov chain on the space of all partitions of the number of cards $n$. Popov proved that this game too converges to triangular configurations, in the sense that the stationary probability measure of the set of configurations that are close to triangular (with a slope that depends on $p$) is close to 1.

There are many other possible ways to formulate stochastic versions of Bulgarian solitaire, but to our knowledge no other possibility has been studied before. Here we consider a particularly natural version, where selection acts on cards rather than piles: When forming a new pile by picking cards from the old piles, every card is picked with a fixed probability $0<p<1$, independently of all other cards.  Our main aim is to show that this \emph{card-based stochastic Bulgarian solitaire} does not converge to a triangular configuration but to an exponential shape. The intuitive reason for this difference in results is that the speed by which a pile loses its cards is constant in the pile-based solitaire but proportional to the current size of the pile in the card-based version. 

Compared to Popov's treatment of the pile-based stochastic Bulgarian solitaire, our treatment of the card-based version has both similarities and novel elements. The card-based rule makes the total number of picked cards in each round independent of the number of piles. This makes the analysis of the resulting stochastic process less complicated than for the pile-based version. In exchange, we can can make a more refined analysis. Specifically, we investigate a question that does not seem to have been studied in the limit-shape literature before: For a fixed number of steps in some stochastic process on integer partitions, how likely is it that this sequence of steps will end up close to the limiting shape of the stationary distribution?

\begin{figure}[h]
\setlength{\unitlength}{0.08cm}
\centering
\begin{picture}(153,70)

% l b r t
\put(0,0){\includegraphics[scale=0.7,clip=true,trim=60 285 40 335]{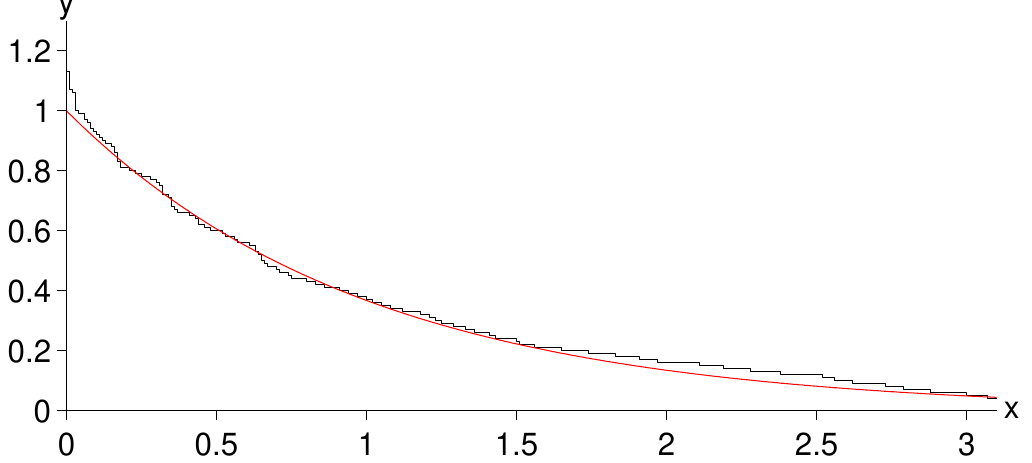}}
\end{picture}
\caption{The result of a computer simulation after 200 moves of the card-based Bulgarian solitaire with $n=10^5$ cards and $p=0.01$ with scaling $1/p$, starting from a triangular configuration. The jagged curve is the diagram-boundary function of the card configuration (defined in Section~\ref{sec:repr_partitions}) and the smooth curve is the limit shape $y=e^{-x}$. \label{fig:sim}}
\end{figure}

The precise result comprises two theorems and a corollary. Theorem~\ref{thm:averageshape} limits the probability that the number of cards in each pile is far from its expected value \emph{after a fixed number of moves}. Theorem~\ref{thm:main} says that as the number of cards tends to infinity the probability of reaching an almost exponential shape \emph{when playing sufficiently long} tends to 1. The outcome of a computer simulation of the process can be seen in Figure~\ref{fig:sim}. In addition, we also consider the \emph{stationary distribution} of the card-based stochastic Bulgarian solitaire regarded as a Markov chain. Corollary~\ref{thm:statmeasure}) says that the stationary probability measure of the set of partitions that have a close to exponential shape is close to 1. 

The remainder of this paper is organized as follows. Section 2.1 introduces notation and makes the solitaire's connection to integer partitions and Young diagrams. Section 2.2 defines what we mean by  a limit shape. Sections 2.3 and 2.4 introduce a key ingredient in our proof, namely the representation of configurations by weak integer compositions, and its connection to the integer partition representation. Section 2.5 discusses the stationary measure of the process defined by the card-based stochastic Bulgarian solitaire. The main result, Theorem~\ref{thm:main}, is formulated in section~3, as well as Corollary~\ref{thm:statmeasure} and its proof. Section 4 makes some general observations needed in the proof of Theorem~\ref{thm:main}. Section 5 contains Theorem~\ref{thm:averageshape}, capturing our limit shape result's relation to the number of moves in the solitaire. Section 6 contains the proof of Theorem~\ref{thm:main}. We end with a discussion in section 7.

\section{Notation and preliminaries}
Let us denote by $\mathscr{B}(n,p)$ the card-based stochastic Bulgarian solitaire with $n$ cards and probability $p$ for a card to be chosen into a new pile.
\subsection{Representation by partitions}
\label{sec:repr_partitions}
Let $\Par(n)=\{\lambda: \lambda\vdash n\}$ be the set of partitions of the integer $n$ into integral parts $\lambda_1\geq\lambda_2\geq\ldots\geq\lambda_\ell>0$, where $\ell=\ell(\lambda)$ is the number of parts of the partition $\lambda$, i.e., $\sum_{i=1}^\ell \lambda_i=n$. We write $\lambda = (\lambda_1,\lambda_2,\dotsc,\lambda_\ell)$. For $i>\ell(\lambda)$ it will be convenient to define $\lambda_i=0$.

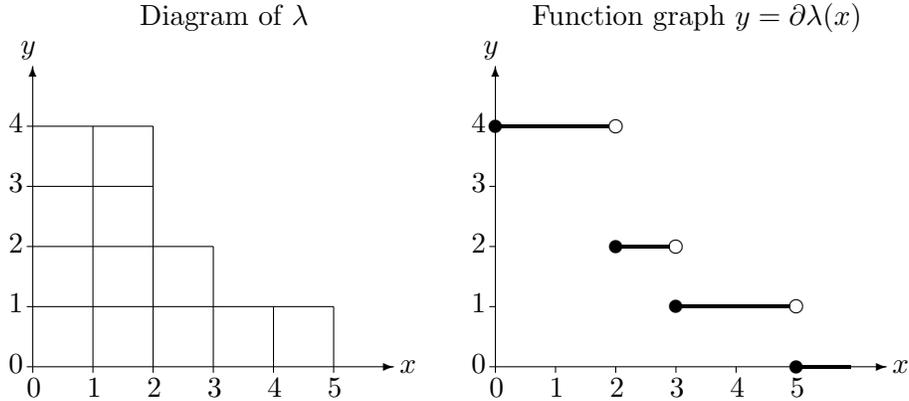
\begin{figure}[h]
\setlength{\unitlength}{0.08cm}
\centering
\begin{tabular}{cc}
Diagram of $\lambda$ & Function graph $y=\partial\lambda(x)$ \\
\begin{picture}(70,60)
\put(4,5){\vector(0,1){50}}
\put(4,5){\vector(1,0){60}}
\put(3,0){$0$}
\put(4,4){\line(0,1){1}}
\put(13,0){$1$}
\put(14,4){\line(0,1){1}}
\put(23,0){$2$}
\put(24,4){\line(0,1){1}}
\put(33,0){$3$}
\put(34,4){\line(0,1){1}}
\put(43,0){$4$}
\put(44,4){\line(0,1){1}}
\put(53,0){$5$}
\put(54,4){\line(0,1){1}}

\put(0,4){$0$}
\put(3,5){\line(1,0){1}}
\put(0,14){$1$}
\put(3,15){\line(1,0){1}}
\put(0,24){$2$}
\put(3,25){\line(1,0){1}}
\put(0,34){$3$}
\put(3,35){\line(1,0){1}}
\put(0,44){$4$}
\put(3,45){\line(1,0){1}}
%\linethickness{1pt}
\put(14,5){\line(0,1){40}}
\put(24,5){\line(0,1){40}}
\put(34,5){\line(0,1){20}}
\put(44,5){\line(0,1){10}}
\put(54,5){\line(0,1){10}}

\put(4,15){\line(1,0){50}}
\put(4,25){\line(1,0){30}}
\put(4,35){\line(1,0){20}}
\put(4,45){\line(1,0){20}}

\put(65,4){$x$}
\put(2,57){$y$}
\end{picture}
\rule{0pt}{60pt} & 
\begin{picture}(75,60)
\put(4,5){\vector(0,1){50}}
\put(4,5){\vector(1,0){65}}
\put(3,0){$0$}
\put(4,4){\line(0,1){1}}
\put(13,0){$1$}
\put(14,4){\line(0,1){1}}
\put(23,0){$2$}
\put(24,4){\line(0,1){1}}
\put(33,0){$3$}
\put(34,4){\line(0,1){1}}
\put(43,0){$4$}
\put(44,4){\line(0,1){1}}
\put(53,0){$5$}

\put(0,4){$0$}
\put(3,5){\line(1,0){1}}
\put(0,14){$1$}
\put(3,15){\line(1,0){1}}
\put(0,24){$2$}
\put(3,25){\line(1,0){1}}
\put(0,34){$3$}
\put(3,35){\line(1,0){1}}
\put(0,44){$4$}

\linethickness{1.3pt}
\put(4,45){\circle*{2}}
\put(24,45){\circle{2}}
\put(4,45){\line(1,0){19}}

\put(24,25){\circle*{2}}
\put(34,25){\circle{2}}
\put(24,25){\line(1,0){9}}

\put(34,15){\circle*{2}}
\put(34,15){\line(1,0){19}}
\put(54,15){\circle{2}}

\put(54,5){\circle*{2}}
\put(54,5){\line(1,0){9}}

\put(70,4){$x$}
\put(2,57){$y$}
\end{picture}
\end{tabular}
\caption{The partition $\lambda=(4,4,2,1,1)\in\Par(12)$.\label{fig:yd_example}}
\end{figure}

We shall represent an integer partition $\lambda$ by a Young diagram drawn as columns of squares in the first quadrant such that the $i$th column has height $\lambda_i$.  For example, the configuration of 12 cards in which there are five piles of sizes 4, 4, 2, 1 and 1 corresponds to the partition $(4,4,2,1,1)\vdash 12$, which is represented by the left diagram in Figure~\ref{fig:yd_example}.

When we speak of shapes of integer partitions we shall mean the shape of the boundary of the Young diagram drawn in this way. To this end, for any partition $\lambda$, define its \textit{diagram-boundary function} as the nonnegative, integer-valued, weakly decreasing and piecewise constant function $\partial\lambda:\Rnn\rightarrow\mathbb{N}$ given by
$$
\partial\lambda(x)=\lambda_{\lfloor x \rfloor+1}.
$$
For example, the right diagram in Figure~\ref{fig:yd_example} depicts the function graph $y=\partial\lambda(x)$ for $\lambda=(4,4,2,1,1)$. (Note that, since we defined $\lambda_i=0$ for $i>\ell(\lambda)$, we have $\partial\lambda(x)=0$ for $x\geq\ell(\lambda)$.) 

As $n$ grows we need to rescale the diagram to achieve any limiting behaviour. Following \cite{Eriksson2012575} and \cite{VershikStatMech}, the diagram is rescaled using some \emph{scaling factor} $a>0$ such that all row lengths are multiplied by $1/a$ and all column heights are multiplied by $a/n$, yielding a constant area of 1.  
Thus, given a partition $\lambda$, define the $a$\textit{-rescaled} diagram-boundary function of $\lambda$ as the positive, real-valued, weakly decreasing and piecewise constant function $\rescaled{a}{\lambda}:\Rnn\rightarrow\Rnn$ given by
\begin{equation}
\label{eq:def_rescaled}
\rescaled{a}{\lambda}(x)=\dfrac{a}{n}\partial\lambda(ax)=\frac{a}{n}\lambda_{\lfloor ax \rfloor+1}.
\end{equation}

\subsection{Limit shapes of birth-and-death processes on Young diagrams}
Eriksson and Sj\"ostrand \cite{Eriksson2012575} studied limit shapes of birth-and-death processes on Young diagrams. In their processes, every step had a single square die and a single square be born. Every move of the card-based stochastic Bulgarian solitaire moves several cards. It can therefore be regarded as a \emph{multi-square} birth-and-death process on Young diagrams. Thus our present study extends the family of processes studied in \cite{Eriksson2012575}.

For each positive integer $n$, let $\nu^{(n)}$ be a
%the stationary probability
distribution on $\Par(n)$.
%obtained from the card-based stochastic Bulgarian solitaire.
(We will study the specific case when $\nu^{(n)}$ is the stationary distribution of the card-based stochastic Bulgarian solitaire. The existence of this distribution will be established in section \ref{subsec:statmeasure}.) We are interested in finding a sequence $\{a_n\}$ of scaling factors such that the rescaled diagrams approach a \textit{limit shape} $\phi$ in probability as $n$ grows to infinity. The precise meaning of this is that, for any $\varepsilon>0$,
\begin{equation}
\label{eq:limitshape_def}
\lim_{n\rightarrow\infty}\nu^{(n)}\left\lbrace \lambda\in\Par(n): \|\rescaled{a_n}{\lambda}(x)-\phi(x)\|_{\infty}<\varepsilon\right\rbrace=1.
\end{equation}
where $\|\cdot\|_{\infty}$ denotes the max-norm $\| f \|_{\infty}=\sup \big\{|f(x)|:x\geq 0\big\}$.
We remark that Vershik \cite{VershikStatMech} and Erlihson and Granovsky \cite{erlihson2008limit} use a weaker condition for convergence towards a limit shape, namely that
$$
\lim_{n\rightarrow\infty}\nu^{(n)}\left\lbrace \lambda\in\Par(n): \sup_{x\in[a,b]} |\rescaled{a_n}{\lambda}(x)-\phi(x)|<\varepsilon\right\rbrace=1
$$
should hold for any compact interval $[a,b]$, and any $\varepsilon>0$. Yakubovich \cite{yakubovich2012ergodicity} and Eriksson and Sjöstrand \cite{Eriksson2012575} use an even weaker condition:
$$
\lim_{n\rightarrow\infty}\nu^{(n)}\left\lbrace \lambda\in\Par(n): |\rescaled{a_n}{\lambda}(x)-\phi(x)|<\varepsilon\right\rbrace=1
$$
for all $x>0$ and all $\varepsilon>0$. 

\subsection{Representation by weak compositions}
\label{sec:repr_compositions}
It will sometimes be convenient to consider a configuration of $n$ cards as a \emph{weak integer composition}, by which we mean an infinite sequence $\alpha=(\alpha_1,\alpha_2,\dotsc)$, not necessarily decreasing, of nonnegative integers adding up to $n$. Let $\mathcal{W}(n)$ be the set of weak compositions of the integer $n$.

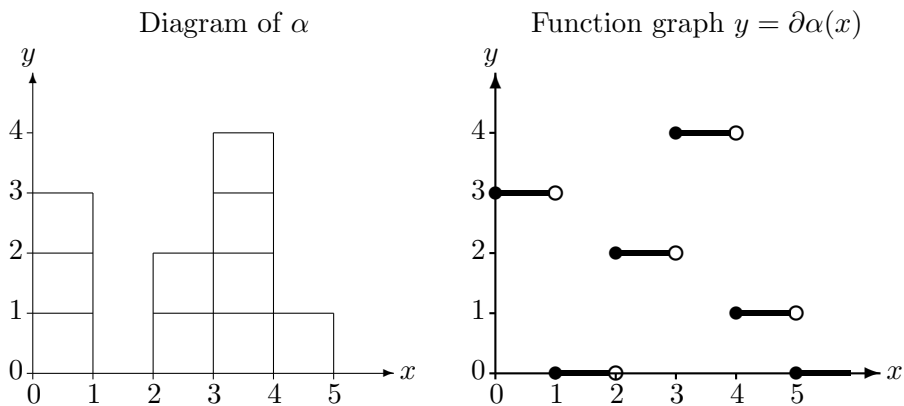
\begin{figure}[h]
\setlength{\unitlength}{0.08cm}
\centering

\begin{tabular}{cc}
Diagram of $\alpha$ & Function graph $y=\partial\alpha(x)$ \\
\begin{picture}(70,60)
%\linethickness{1pt}
\put(4,5){\vector(0,1){50}}
\put(4,5){\vector(1,0){60}}
\put(3,0){$0$}
\put(4,4){\line(0,1){1}}
\put(13,0){$1$}
\put(14,4){\line(0,1){1}}
\put(23,0){$2$}
\put(24,4){\line(0,1){1}}
\put(33,0){$3$}
\put(34,4){\line(0,1){1}}
\put(43,0){$4$}
\put(44,4){\line(0,1){1}}
\put(53,0){$5$}
\put(54,4){\line(0,1){1}}

\put(0,4){$0$}
\put(3,5){\line(1,0){1}}
\put(0,14){$1$}
\put(3,15){\line(1,0){1}}
\put(0,24){$2$}
\put(3,25){\line(1,0){1}}
\put(0,34){$3$}
\put(3,35){\line(1,0){1}}
\put(0,44){$4$}
\put(3,45){\line(1,0){1}}
%\linethickness{1pt}
\put(14,5){\line(0,1){30}}
\put(24,5){\line(0,1){20}}
\put(34,5){\line(0,1){40}}
\put(44,5){\line(0,1){40}}
\put(54,5){\line(0,1){10}}

\put(4,15){\line(1,0){10}}
\put(24,15){\line(1,0){30}}
\put(4,25){\line(1,0){10}}
\put(24,25){\line(1,0){20}}
\put(34,35){\line(1,0){10}}
\put(34,45){\line(1,0){10}}
\put(4,35){\line(1,0){10}}

\put(65,4){$x$}
\put(2,57){$y$}
\end{picture}
\rule{0pt}{60pt} & 
\begin{picture}(75,60)
\thicklines
%\linethickness{1pt}
\put(4,5){\vector(0,1){50}}
\put(4,5){\vector(1,0){64}}

\put(3,0){$0$}
\put(4,4){\line(0,1){1}}
\put(13,0){$1$}
\put(14,4){\line(0,1){1}}
\put(23,0){$2$}
\put(24,4){\line(0,1){1}}
\put(33,0){$3$}
\put(34,4){\line(0,1){1}}
\put(43,0){$4$}
\put(44,4){\line(0,1){1}}
\put(53,0){$5$}

\put(0,4){$0$}
\put(3,5){\line(1,0){1}}
\put(0,14){$1$}
\put(3,15){\line(1,0){1}}
\put(0,24){$2$}
\put(3,25){\line(1,0){1}}
\put(0,34){$3$}
\put(0,44){$4$}
\put(3,45){\line(1,0){1}}

\linethickness{2pt}
\put(4,35){\circle*{2}}
\put(14,35){\circle{2}}
\put(4,35){\line(1,0){9}}

\put(14,5){\circle*{2}}
\put(24,5){\circle{2}}
\put(14,5){\line(1,0){9}}

\put(34,45){\circle*{2}}
\put(44,45){\circle{2}}
\put(34,45){\line(1,0){9}}

\put(24,25){\circle*{2}}
\put(34,25){\circle{2}}
\put(24,25){\line(1,0){9}}

\put(44,15){\circle*{2}}
\put(54,15){\circle{2}}
\put(44,15){\line(1,0){9}}

\put(54,5){\circle*{2}}
\put(54,5){\line(1,0){9}}

\put(69,4){$x$}
\put(2,57){$y$}
\end{picture}
\end{tabular}
\caption{The composition $\alpha=(3,0,2,4,1,0,0,\dotsc)\in\mathcal{W}(10)$.\label{fig:psi_alpha_example}}
\end{figure}

We define the \textit{diagram}, the \textit{boundary function} $\partial\alpha$, and the \textit{rescaled boundary function} $\rescaled{a}{\alpha}$ of a weak integer composition $\alpha$ in exact analogy to the way we defined them for integer partitions in Section~\ref{sec:repr_partitions}. For example,  the diagram of $\alpha=(3,0,2,4,1,0,0,\dotsc)$ and the corresponding function graph $y=\partial\alpha(x)$ are shown in Figure~\ref{fig:psi_alpha_example}.

Moreover, for a stochastic process with state-space $\mathcal{W}(n)$ having a unique stationary distribution $\pi$, we define the \emph{average shape} of this process as the piecewise constant function $\partial_E:\mathbb{R}_{>0} \rightarrow \mathbb{R}_{\ge 0}$ given by
$$
\partial_E(x) = E\alpha_{\lfloor x \rfloor + 1},
$$
where $\alpha\in\mathcal{W}(n)$ is sampled from $\pi$.

\subsection{A necessary lemma}
\label{sec:thelemma}
In the proof of the limit shape result, Theorem~\ref{thm:main}, we need to be able to keep track of individual piles of cards. We shall therefore order piles by time of creation rather than by size, which means that configurations are represented by weak integer compositions rather than integer partitions.  It will turn out that the rescaled diagram-boundary function of these compositions tends to a limit function as the game is played and as the number of cards
tends to infinity. However, in the end we want to express the limit-shape result in terms of diagram-boundary functions of integer partitions, not weak integer compositions. For any $\alpha\in\mathcal{W}(n)$, define the operator $\ord$ as the ordering operator that arranges the parts of $\alpha$ in descending order. If we omit trailing zeros, $\ord\alpha$ is an integer partition of $n$. For example, if $\alpha=(3,0,2,4,2,0,0,\dotsc)$, then $\ord\alpha=(4,3,2,2)$. In this section, we prove an important lemma which says that such sorting of the piles by size does not harm the convergence to a limiting shape.

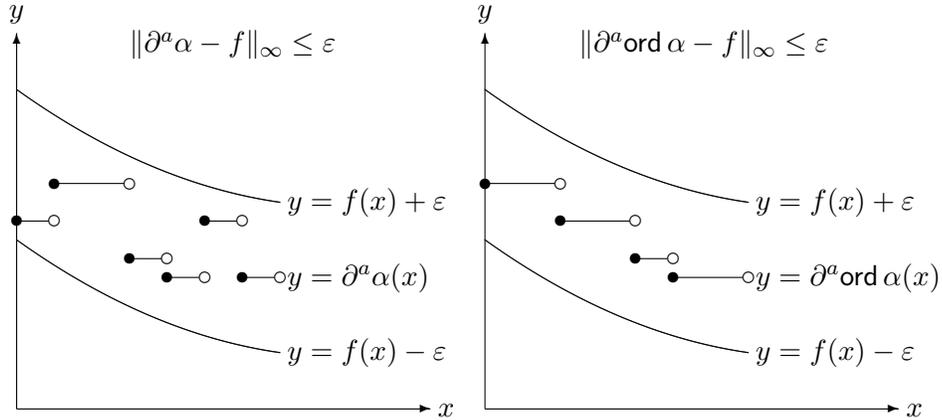
\begin{figure}[h]
\setlength{\unitlength}{0.05cm}
\centering
\begin{tabular}{cc}
\begin{picture}(110,105)
% Coordinate system
\put(0,0){\vector(1,0){110}}
\put(0,0){\vector(0,1){100}}
\put(112,-2){$x$}
\put(-2,104){$y$}
% Curves
\qbezier(0,85)(35,60)(70,55)
\qbezier(0,45)(35,20)(70,15)
% Equations
\put(30,95){$\| \rescaled{a}{\alpha}-f \|_{\infty} \leq \varepsilon$}
\put(72,53){$y=f(x)+\varepsilon$}
\put(72,33){$y=\rescaled{a}{\alpha}(x)$}
\put(72,13){$y=f(x)-\varepsilon$}
% Steps
\put(1,50){\line(1,0){7.5}}
\put(0,50){\circle*{3}}
\put(10,50){\circle{3}}

\put(11,60){\line(1,0){17.5}}
\put(10,60){\circle*{3}}
\put(30,60){\circle{3}}

\put(31,40){\line(1,0){7.5}}
\put(30,40){\circle*{3}}
\put(40,40){\circle{3}}

\put(41,35){\line(1,0){7.5}}
\put(40,35){\circle*{3}}
\put(50,35){\circle{3}}

\put(51,50){\line(1,0){7.5}}
\put(50,50){\circle*{3}}
\put(60,50){\circle{3}}

\put(61,35){\line(1,0){7.5}}
\put(60,35){\circle*{3}}
\put(70,35){\circle{3}}

\end{picture}
\hspace{5pt} &
%------------------------------------------------------------
\begin{picture}(110,105)
% Coordinate system
\put(0,0){\vector(1,0){110}}
\put(0,0){\vector(0,1){100}}
\put(112,-2){$x$}
\put(-2,104){$y$}
% Curves
\qbezier(0,85)(35,60)(70,55)
\qbezier(0,45)(35,20)(70,15)
% Equations
\put(25,95){$\| \rescaled{a}{\ord\alpha}-f \|_{\infty} \leq \varepsilon$}
\put(72,53){$y=f(x)+\varepsilon$}
\put(72,33){$y=\rescaled{a}{\ord\alpha}(x)$}
\put(72,13){$y=f(x)-\varepsilon$}
% Steps
\put(1,60){\line(1,0){17.5}}
\put(0,60){\circle*{3}}
\put(20,60){\circle{3}}

\put(21,50){\line(1,0){17.5}}
\put(20,50){\circle*{3}}
\put(40,50){\circle{3}}

\put(41,40){\line(1,0){7.5}}
\put(40,40){\circle*{3}}
\put(50,40){\circle{3}}

\put(51,35){\line(1,0){17.5}}
\put(50,35){\circle*{3}}
\put(70,35){\circle{3}}

\end{picture}
\end{tabular}
\caption{An example of a composition $\alpha$ and a decreasing function $f$ showing that if $\rescaled{a}{\alpha}(x)$ is enclosed between $f(x)-\varepsilon$ and $f(x)+\varepsilon$, then so is $\rescaled{a}{\ord\alpha}(x)$, an immediate consequence of Lemma~\ref{lem:reordering}.  \label{fig:reordering_proof}}
\end{figure}

The proof of the lemma will use some basic theory of symmetric-decrea\-sing rearrangements, see for example  \cite[Ch.~10]{inequalities}
or \cite[Ch.~3]{analysis}. For any measurable function $f\colon\mathbb{R}\rightarrow\Rnn$ such that $\lim_{x\rightarrow\pm\infty}f(x)=0$, there is a
unique function $f^\ast\colon\mathbb{R}\rightarrow\Rnn$, called the \emph{symmetric-decreasing rearrangement} of $f$,
with the following properties:
\begin{itemize}
\item $f^\ast$ is symmetric, that is, $f^\ast(-x)=f^\ast(x)$ for all $x$,
\item $f^\ast$ is weakly decreasing on the interval $[0,\infty)$,
\item $f^{\ast}$ and $f$ are equimeasurable, that is,
$$\mathcal{L}(\{ x:\;f(x)>t \})=\mathcal{L}(\{ x:\;f^\ast(x)>t \})$$ for all $t>0$, where $\mathcal{L}$ denotes the Lebesgue measure,
\item $f^{\ast}$ is lower semicontinuous.
\end{itemize}
In particular, if $f$ is a symmetric function that is weakly decreasing and right-continuous on $[0,\infty)$ and tends to $0$ at infinity, then $f^\ast=f$.
\begin{lemma}
\label{lem:reordering}
Let $\alpha\in\mathcal{W}(n)$ be a weak composition of $n$, let $a>0$ be any scaling factor and let $f:\Rnn\rightarrow \Rnn$ be a right-continuous and weakly decreasing function such that $f(x)\rightarrow 0$ as $x\rightarrow\infty$. The rescaled diagram-boundary functions before and after sorting of the weak composition satisfy the inequality
$$
\| \rescaled{a}{\ord \alpha}-f\|_{\infty}\le
\| \rescaled{a}{\alpha}-f\|_{\infty}.
$$
%$\psi_{\textup{ord}(\alpha)}\in\mathcal{N}_f(\varepsilon,n)$.
\end{lemma}
\begin{proof}
The intuition of the lemma should be obvious from Figure~\ref{fig:reordering_proof}. To be able to use the standard machinery of symmetric rearrangements, we consider the functions
$f$, $\rescaled{a}{\alpha}$, and $\rescaled{a}{\ord \alpha}$ as being defined
on the entire real axis by letting $f(x)=f(|x|)$ and analogously for $\rescaled{a}{\alpha}$,
and $\rescaled{a}{\ord \alpha}$.

Since $f(x)\rightarrow 0$ as $x\rightarrow\infty$, its symmetric-decreasing rearrangement $f^{\ast}$ is defined, and,
since $f$ is weakly decreasing and lower semicontinuous, we have $f^{\ast}=f$.
Similarly, $\rescaled{a}{\ord\alpha}(x)\rightarrow 0$ as $x\rightarrow\infty$ and is weakly decreasing, so
$(\rescaled{a}{\ord\alpha})^{\ast}=\rescaled{a}{\ord\alpha}$.
Moreover,  $(\rescaled{a}{\alpha})^{\ast}=\rescaled{a}{\ord\alpha}$ must hold because the operator $\ord$ arranges the composition
parts in descending order.

Now, since symmetric rearrangements decrease $L^p$-distances for any $1\le p\le\infty$ (see for example \cite{analysis}, Section 3.4), we obtain
\[
\| \rescaled{a}{\ord \alpha}-f\|_{\infty}=\| (\rescaled{a}{\alpha})^{\ast}-f^{\ast}\|_{\infty}\le
\| \rescaled{a}{\alpha}-f\|_{\infty}.
\]
\end{proof}

\begin{corollary}
For any distribution $\rho^{(n)}$ on $\mathcal{W}(n)$, define a corresponding distribution $\tilde{\rho}^{(n)}$ on $\Par(n)$  by
\begin{equation}
\label{eq:reordering_cor}
\tilde{\rho}^{(n)}(\lambda)=\sum_{\substack{ \alpha\in\mathcal{W}(n) \\ \ord\alpha=\lambda }}\rho^{(n)}(\alpha).
\end{equation}
If $\phi$ is a limit shape of $\rho^{(n)}$ on $\mathcal{W}(n)$ with some scaling $\{a_n\}$, then $\phi $ is also a limit shape of $\tilde{\rho}^{(n)}$ on $\Par(n)$ with the same scaling.
\end{corollary}
%Formally, $\nu^{(n)}$ depends on $\rho^{(n)}$, but for simplicity of notation, we do not indicate in $\nu^{(n)}$ the dependence on $\rho^{(n)}$, but assume the underlying distribution on $\mathcal{W}(n)$. 
\begin{proof}
The assumption that $\phi$ is a limit shape of the distribution $\rho^{(n)}$ on $\mathcal{W}(n)$ with scaling $\{a_n\}$ means that
$$
\lim_{n\rightarrow\infty}\rho^{(n)}\left\lbrace \alpha\in\mathcal{W}(n): \|\rescaled{a_n}{\alpha}-\phi\|_{\infty}<\varepsilon\right\rbrace=1.
$$
By virtue of Lemma~\ref{lem:reordering} we can replace $\alpha$ with $\ord\alpha$ in this formula:
\begin{equation}
\label{eq:lemmaproof1}
\lim_{n\rightarrow\infty}\rho^{(n)}\left\lbrace \alpha\in\mathcal{W}(n): \|\rescaled{a_n}{\ord\alpha}-\phi\|_{\infty}<\varepsilon\right\rbrace=1.
\end{equation}
The set of weak compositions $A:=\left\lbrace \alpha\in\mathcal{W}(n): \|\rescaled{a_n}{\ord\alpha}-\phi\|_{\infty}<\varepsilon\right\rbrace$ can be written as a disjoint union of equivalence classes with respect to sorting:
$$
A=\bigcup_{\lambda\in L} \{ \alpha\in\mathcal{W}(n):\, \ord\alpha=\lambda \}.
$$
where $L=\{ \lambda\in\Par(n):\, \|\rescaled{a_n}{\lambda}-\phi\|_{\infty}<\varepsilon \}$. The $\rho^{(n)}$-probability measure of $A$ is
\begin{align*}
\rho^{(n)}(A)
& = \rho^{(n)}\left( \bigcup_{\lambda\in L} \{ \alpha\in\mathcal{W}(n):\, \ord\alpha=\lambda \} \right) \\
& = \sum_{\lambda\in L} \rho^{(n)}\{ \alpha\in\mathcal{W}(n):\, \ord\alpha=\lambda \} \\
& = \sum_{\lambda\in L} \; \sum_{\substack{  \alpha\in\mathcal{W}(n) \\ \ord\alpha=\lambda }} \rho^{(n)}(\alpha) \\
& = \sum_{\lambda\in L} \tilde{\rho}^{(n)}(\lambda) && \text{(by \eqref{eq:reordering_cor})} \\
& = \tilde{\rho}^{(n)}(L).
\end{align*}
From \eqref{eq:lemmaproof1} we have that $\lim_{n\rightarrow\infty}\rho^{(n)}(A)=1$. Because $\rho^{(n)}(A)=\tilde{\rho}^{(n)}(L)$, we can conclude that also ${\lim_{n\rightarrow\infty}\tilde{\rho}^{(n)}(L)=1}$, that is,
$$
\lim_{n\rightarrow\infty}\tilde{\rho}^{(n)}\left\lbrace \lambda\in\mathcal\Par(n): \|\rescaled{a_n}{\lambda}-\phi\|_{\infty}<\varepsilon\right\rbrace=1.
$$
This means that $\phi$ is a limit shape of the distribution $\nu^{(n)}$ on $\Par(n)$.
\end{proof}

\subsection{The stationary measure}
\label{subsec:statmeasure}
The card-based stochastic Bulgarian solitaire can be regarded as a Markov chain with state-space $\Par{(n)}$. Let us denote it by $(\lambda^{(0)},\lambda^{(1)},\dotsc)$. This Markov chain is aperiodic and irreducible. It is irreducible because, starting from any state, an arbitrary $\lambda\in\mathcal{P}(n)$ can be reached in $\ell(\lambda)$ moves as follows. In the first move, choose  cards for the last column in $\lambda$; in the next move, choose (among remaining cards) the cards for the next-to-the-last column in $\lambda$; repeat until all the $\ell(\lambda)$ columns have been chosen. All these selections of cards have probability $>0$, hence the irreducibility. It is aperiodic because all states are aperiodic: There is always a positive probability to choose zero cards in a move and hence remain in the same state. 

It is well known that a finite state-space irreducible Markov chain has a unique stationary distribution $\pi$, and that, if it also is aperiodic, the distribution converges to $\pi$ starting from any initial state. We denote by $\pi_{p,n}$ the stationary measure of the Markov chain $(\lambda^{(0)},\lambda^{(1)},\dotsc)$ on $\Par(n)$  given by the card-based stochastic Bulgarian solitaire.

Readers acquainted with the limit shape literature may wonder whether the stationary measure $\pi_{p_n,n}$ has the property of being \emph{multiplicative}, in the sense of interpretable as the product measure on the space of integer sequences restricted to a certain affine subspace \cite{fristedt1993structure}. The multiplicative property is useful in limit shape problems and related problems  \cite{VershikStatMech,corteel1999multiplicity, erlihson2004reversible, goh2008random, pittel1997likely}. However,  such techniques will not be used here as $\pi_{p_n,n}$ is unlikely to be multiplicative in general. 

\section{The limit shape result}
We shall now state the limit shape result. If we view the solitaire as a process on Young diagrams, Theorem~\ref{thm:main} says that after a sufficiently large number $m$ of moves the $(1/p)$-rescaled boundary function of the diagram will resemble the exponential shape $e^{-x}$ with high probability. 

\begin{theorem}
\label{thm:main}
For each positive integer $n$, pick a probability $p_n\in(0,1)$, a (possibly random)
initial configuration $\lambda^{(0)}\in\Par(n)$ and let $(\lambda^{(0)},\lambda^{(1)},\dotsc)$ be the Markov chain on $\Par(n)$ defined by $\mathscr{B}(n,p_n)$.

Suppose that $p_n$ as a function of $n$ has the asymptotical properties that
\begin{itemize}
\item
$p_n\rightarrow 0$ as $n\rightarrow\infty$ and
\item
$p_n=\omega(\frac{\log n}{n})$, i.e. $p_n\big/\frac{\log n}{n} \rightarrow \infty$ as $n\rightarrow\infty$.
\end{itemize}

Then, for any 
$m=m(\varepsilon_1,n)>\frac{\varepsilon^2}{2+\varepsilon}n$,
the probability distribution for the resulting diagram $\lambda^{(m)}$ after playing $m$ moves has the limit shape $g(x)=e^{-x}$ with scaling $1/p_n$.
In fact,
for any $\varepsilon>0$ we have
$$
P\left(\| \rescaled{1/p_n}{\lambda^{(m)}}-g\|_{\infty} \leq \varepsilon \right)
\geq 1-\exp\left[-\frac{\varepsilon^2}{2+\varepsilon}np_n\bigl(1-o(1)\bigr)\right].
$$
\end{theorem}

As a simple consequence of Theorem~\ref{thm:main}, the stationary distribution of the card-based Bulgarian solitaire also has the limit shape $e^{-x}$.
\begin{corollary}
\label{thm:statmeasure}
Let $p_n$ be a probability, dependent on $n$, with the same asymptotical properties as in Theorem~\ref{thm:main} and let $\pi_{p_n,n}$ denote the stationary measure of the Markov chain $(\lambda^{(0)},\lambda^{(1)},\dotsc)$ on $\Par(n)$ defined by $\mathscr{B}(n,p_n)$. Then $\pi_{p_n,n}$ has the limit shape $g(x)=e^{-x}$ with scaling $1/p_n$. In fact, for any $\varepsilon>0$ we have
$$
\pi_{p_n,n}\bigl(\{ \lambda\in\Par(n): \| \rescaled{1/p_n}{\lambda}-g \|_{\infty} \leq \varepsilon \})\geq 1-\exp\left[-\frac{\varepsilon^2}{2+\varepsilon}np_n\bigl(1-o(1)\bigr)\right].
$$
\end{corollary}
\begin{proof}
Since $\pi_{p_n,n}$ is the stationary distribution, if we start with a partition $\lambda^{(0)}$ sampled from $\pi_{p_n,n}$ and play $m$ moves, the resulting partition $\lambda^{(m)}$ will also be sampled from $\pi_{p_n,n}$. Thus, the corollary follows from Theorem~\ref{thm:main} by choosing $\lambda^{(0)}$ as a stochastic partition sampled from the stationary distribution.
\end{proof}
Let us discuss what happens if $p_n$ does not fulfil the conditions in the theorem. If $p_n$ is bounded away from zero, then the scaling $1/p_n$ is bounded and hence cannot transform the jumpy boundary diagrams into a continuous limit shape. On the other hand,
if $p_n$ tends to zero too fast, the pile sizes will be small and their random fluctuations will be large. For instance, the new pile after each move has a size
drawn from the binomial distribution $\Bin(n,p_n)$ with relative standard deviation
$\sim 1/\sqrt{np_n}$.

\section{Some observations}
The proof of Theorem~\ref{thm:main} will rely on some general analytic and probabilistic observations that we present in this section.

The first observation is that the well-known convergence of $\left( 1-\frac{1}{n} \right)^{nx}$ to $e^{-x}$ as $n\rightarrow\infty$ is uniform on $x\geq 0$. Setting $p=\frac{1}{n}$, this result can be states as followis. 
\begin{observation}
\label{prop:uniformconv}
The convergences
\begin{align}
\label{eq:uniformconv1}
\left( 1-p \right)^{x/p} - e^{-x} & \rightarrow 0 \quad\text{and}\quad \\
\label{eq:uniformconv2}
\left( 1-p \right)^x - e^{-px} & \rightarrow 0 
\end{align}
hold uniformly on the interval $x\in[0,\infty)$ as $p\rightarrow 0$.
\end{observation}

\begin{proof}
Let $f(p)=(1-p)^{1/p}$. Since
the Taylor expansion at $p=0$ of $\log f(p)=\frac{1}{p}\log(1-p)$ is
\begin{equation*}
\log f(p)=-1-\frac{p}{2}-\frac{p^2}{3}-\dotsb,
\end{equation*}
the function $-\log f(p)$ tends to 1 from above as $p\searrow 0$, and hence for any $y\in[0,1]$, we have $y^{-\log f(p)}\nearrow y$ as $p\searrow 0$. Since $[0,1]$ is a compact set, Dini's Theorem can be applied. Thus  $y^{-\log f(p)}-y\rightarrow 0$ uniformly on $y\in[0,1]$ as $p\rightarrow 0$. Substituting $y=e^{-x}$ yields \eqref{eq:uniformconv1} while instead substituting $y=e^{-px}$ yields \eqref{eq:uniformconv2}.
\end{proof}

We shall also need the following version of Chernoff bounds.

\begin{observation}
\label{cor:chernoff_abs}
For $n\ge 0$ and $0\le p\le 1$, let $X\sim\Bin(n,p)$ and set $\mu=E(X)=np$. Then, for any $\gamma\geq 0$,
\begin{equation*}
%\label{eq:chernoff_abs}
P(\left|X-\mu\right|\geq \gamma) \leq
2\exp\left(-\frac{\gamma^2}{2\mu+\gamma}\right).
\end{equation*}
\end{observation}
\begin{proof}
By setting $\gamma=\eta\mu$ for $\eta\ge 0$, we can rewrite this bound in terms of the relative deviation $\eta$ instead of the absolute $\gamma$:
\begin{equation}
\label{eq:chernoffproof1}
P(\left|X-\mu\right|\geq \gamma) = P(X\le (1-\eta)\mu) + P(X\ge(1+\eta)\mu)
\end{equation}
Now we can use known bounds on each of the terms in the right hand side of \eqref{eq:chernoffproof1}. By Theorem 4.5 in \cite{mitzenmacher2005probability}, the first term has the bound
\begin{equation}
\label{eq:chernoffproof_term1}
P(X\leq(1-\eta)\mu) \leq e^{-\frac{\eta^2}{2}\mu} = \exp\left(-\frac{\gamma^2}{2\mu}\right).
\end{equation}
By Theorem 4.4 in \cite{mitzenmacher2005probability}, the second term has the bound
\begin{equation*}
%\label{eq:chernoffproof_term2}
P(X\geq(1+\eta)\mu) \le \left( \frac{e^{\eta}}{(1+\eta)^{1+\eta}} \right)^{\mu} = 
\left( \frac{e^{\eta}}{e^{(1+\eta)\log(1+\eta)}} \right)^{\mu}.
\end{equation*}
Using the inequality $\log(1+\eta) \ge \frac{2\eta}{2+\eta}$ for $\eta \ge 0$ in the right hand side, we get
%of \eqref{eq:chernoffproof1}
\begin{equation}
\label{eq:chernoffproof_term2}
P(X\geq(1+\eta)\mu) \le \left( \frac{e^{\eta}}{e^{(1+\eta)\frac{2\eta}{2+\eta}}} \right)^{\mu} = 
e^{-\frac{\eta^2}{2+\eta}\mu}=
\exp\left(-\frac{\gamma^2}{2\mu+\gamma}\right).
\end{equation}
Using the inequalities \eqref{eq:chernoffproof_term1} and \eqref{eq:chernoffproof_term2} and the fact that $\exp\left(-\frac{\gamma^2}{2\mu}\right) \leq \exp\left(-\frac{\gamma^2}{2\mu+\gamma}\right)$ for $\gamma\ge 0$, we can write \eqref{eq:chernoffproof1} as
$$
P(\left|X-\mu\right|\geq \gamma) \le 2\exp\left(-\frac{\gamma^2}{2\mu+\gamma}\right),
$$
which is what we wanted to prove.
\end{proof}

\section{The average shape result}
We will now state Theorem~\ref{thm:averageshape} which, loosely speaking, limits the probability that after a \emph{fixed} number of moves $m$, the pile-sizes in the card-based Bulgarian solitaire are far from their expected value. More precisely, when representing a card configuration of the solitaire as a weak integer composition $\alpha$, after $m$ moves in the solitaire, the boundary function $\partial\alpha$ resembles the average shape $\partial_E$ with probability almost 1. We explicitly relate the error in probability to the number of moves $m$.

In Corollary~\ref{cor:averageshape}, we will then choose $m$ in such a way that we can use it in the proof of the limit shape result, Theorem~\ref{thm:main}.
\begin{theorem}
\label{thm:averageshape}
Let $(\alpha^{(0)},\alpha^{(1)},\dotsc)$ be the Markov chain on $\mathcal{W}(n)$ defined by $\mathscr{B}(n,p)$ with a (possibly random) initial configuration $\alpha^{(0)}\in\mathcal{W}(n)$. Then, for any $\varepsilon_1>0$ and any number of moves $m$, we have
\begin{equation}
\label{eq:avshaperesult1}
P\bigl( \forall k \leq m : |\alpha_k^{(m)}-E\alpha_k^{(m)}| < \varepsilon_1 np \bigr) > 1-2m\exp\left( - \frac{\varepsilon_1^2}{2+\varepsilon_1}np \right)
\end{equation}
and
\begin{equation}
\label{eq:avshaperesult2}
P\bigl(\forall k>m: \alpha_k^{(m)}=0\bigr) \geq 1-n(1-p)^m.
\end{equation}
\end{theorem}
\begin{proof}
In order to keep track of the piles,  put each pile in a bowl and line up the bowls in a row on the table. In each move of the game, the new (possibly empty) pile is put in a new bowl to the left of all old bowls. Since the theorem is only concerned with the configuration $\alpha^{(m)}$ after $m$ moves, let us omit the superscript $(m)$. Let $\alpha_k$ be the number of cards in the $k$th bowl \emph{from the left} after $m$ moves (i.e., $\alpha_1$ is always the number of cards in the most recently formed pile). For convenience, let there be infinitely many empty bowls to the right. This way, a configuration of cards corresponds to a weak integer composition $\alpha=(\alpha_1,\alpha_2,\dotsc)$, as described in section~\ref{sec:repr_compositions}.

We shall consider two separate regimes of $k$-values:  $k\le m$  and $k>m$.

\underline{Regime 1}: $k\leq m$. In order for a specific card to be in the $k$th bowl (i) it must have been picked in the move when the bowl was created, for which the probability is $p$, and (ii) it must have stayed in that bowl and not been picked in the following $k-1$ moves, for which the probability is $(1-p)^{k-1}$. Thus,
\begin{equation}
\label{eq:Ebetak}
\alpha_k\sim\Bin\left(n,p(1-p)^{k-1}\right) \quad\text{and hence}\quad E\alpha_k=np(1-p)^{k-1}.
\end{equation}
Since $\alpha_k$ is a binomially distributed variable, we can use Observation~\ref{cor:chernoff_abs} with $X=\alpha_k$ and $\gamma=\varepsilon_1 np$ to obtain
\begin{align}
P\left(\left|\alpha_k-E\alpha_k\right|\geq\varepsilon_1 np\right) & \leq 
2\exp\left( -\frac{\varepsilon_1^2}{2(1-p)^{k-1}+\varepsilon_1}np \right) & \nonumber \\
& \leq 2\exp\left( -\frac{\varepsilon_1^2}{2+\varepsilon_1}np \right). \label{eq:regime1_onek}
\end{align}
We shall now bound the probability that $|\alpha_k-E\alpha_k|\geq \varepsilon_1 np$ for at least one $k$ in \textit{the entire} first regime (i.e., for at least one value of $k\leq m
%=\lfloor nf(\varepsilon_1) \rfloor$
$) by summing $m$ terms of the type in the right hand side of \eqref{eq:regime1_onek}. This is particularly easy to do as these terms are independent of $k$.
\begin{align*}
& P\left( \exists k \leq m : \left|\alpha_k-E\alpha_k\right| \geq\varepsilon_1 np \right) \\
& \hspace{20pt} \leq \sum_{k=1}^{m} P\left(\left|\alpha_k-E\alpha_k\right| \geq\varepsilon_1 np\right) \\
& \hspace{20pt} \leq m\cdot 2\exp\left( -\frac{\varepsilon_1^2}{2+\varepsilon_1}np \right) && \text{(by }\eqref{eq:regime1_onek}\text{)}
\end{align*}
Therefore, for the complementary event, we have
\begin{equation}
\label{eq:regime1end}
P\bigl( \forall k \leq m : \left|\alpha_k-E\alpha_k\right| < \varepsilon_1 np \bigr) > 1-2m\exp\left( -\frac{\varepsilon_1^2}{2+\varepsilon_1}np \right),
\end{equation}
which proves \eqref{eq:avshaperesult1}.

\underline{Regime 2}: $k>m$. If all $n$ cards have been picked at least once after $m$ moves, all bowls in this regime are empty, i.e.\ $\alpha_k=0$ for all $k>m$. The probability for this is $(1-(1-p)^m)^n$, since the probability that a specific card has not been picked after $m$ moves is $(1-p)^m$. Thus,
$$
P(\forall k>m: \alpha_k=0) = \bigl(1-(1-p)^m\bigr)^n \geq 1-n(1-p)^m,
$$
which proves \eqref{eq:avshaperesult2}.
\end{proof}
It is possible to choose $\varepsilon_1$ in Theorem~\ref{thm:averageshape} in such a way that the errors in probability in \eqref{eq:avshaperesult1} and \eqref{eq:avshaperesult2} are balanced in the sense that they will have the same asymptotic behaviour as $n\rightarrow\infty$. We  achieve this in Corollary~\ref{cor:averageshape}.

\begin{corollary}
\label{cor:averageshape}
Let $p_n$ be a probability, dependent on $n$, with the same asymptotical properties as in Theorem~\ref{thm:main}.
Let $(\alpha^{(0)},\alpha^{(1)},\dotsc)$ be the Markov chain on $\mathcal{W}(n)$ defined by $\mathscr{B}(n,p)$ with a (possibly random) initial configuration $\alpha^{(0)}\in\mathcal{W}(n)$. Let $\varepsilon_1>0$ and define the function $f(x)=\frac{x^2}{2+x}$. Then, for any $m=m(\varepsilon_1,n)>f(\varepsilon_1)n$, we have
\begin{equation}
\label{eq:avshaperesult1cor}
P\bigl( \forall k \leq m : |\alpha_k^{(m)}-E\alpha_k^{(m)}| < \varepsilon_1 np \bigr) > 1-\exp\left[-f(\varepsilon_1)np(1-o(1)) \right]
\end{equation}
and
\begin{equation}
\label{eq:avshaperesult2cor}
P\bigl(\forall k>m: \alpha_k^{(m)}=0\bigr) > 1-\exp\left[-f(\varepsilon_1)np(1-o(1)) \right].
\end{equation}
where $o(1)$ is with respect to $n\rightarrow\infty$.
\end{corollary}
\begin{proof}
First, note that  $nf(\varepsilon_1)<nf(\varepsilon)$. For $m>\lfloor nf(\varepsilon)\rfloor$, 
playing $m$ moves of the solitaire from the initial state $\alpha^{(0)}$ is equivalent to playing $\lfloor nf(\varepsilon_1)\rfloor$ moves from the initial state $\alpha^{(m-\lfloor nf(\varepsilon_1)\rfloor)}$. Without loss of generality, we can therefore assume $m=\lfloor nf(\varepsilon_1)\rfloor$. Also, in the following let us abbreviate $p=p_n$.

Now, let us investigate the (complement to the) probability in \eqref{eq:avshaperesult1cor}:
\begin{align*}
& P\left( \exists k \leq m : |\alpha_k^{(m)}-E\alpha_k^{(m)}| \geq\varepsilon_1 np \right) \\
& \hspace{10pt} \leq 2m\exp\left( -f(\varepsilon_1)np \right) & \text{(by Theorem~\ref{thm:averageshape})} \\
& \hspace{10pt} \leq 2nf(\varepsilon_1)\exp\left( -f(\varepsilon_1)np \right) & \text{(since }m=\lfloor nf(\varepsilon_1)\rfloor<nf(\varepsilon_1)\text{)}\\
& \hspace{10pt} = \exp\left[-f(\varepsilon_1)np\left(1-\frac{\log\bigl(2nf(\varepsilon_1)\bigr)}{f(\varepsilon_1)np} \right) \right] \\
& \hspace{10pt} = \exp\biggl[ -f(\varepsilon_1)np\bigl(1-o(1)\bigr) \biggr] & \text{(since } p=\omega\left(\dfrac{\log n}{n}\right) \text{)},
%& \hspace{10pt} = P_{\varepsilon_1},
\end{align*}
which proves \eqref{eq:avshaperesult1cor}.

Next, we turn our attention to the probability in \eqref{eq:avshaperesult2cor}:
\begin{align}
& P(\forall k>m: \alpha_k^{(m)}=0) \geq 1-n(1-p)^m & \text{(by Theorem \ref{thm:averageshape})} \nonumber \\
& \hspace{8pt} > 1-ne^{-mp} & \text{(since }e^{-p} > 1-p\text{)} \nonumber \\
& \hspace{8pt} = 1-n\exp\bigl(-\lfloor nf(\varepsilon_1)\rfloor p\bigr) \nonumber \\
& \hspace{8pt} = 1-\exp\biggl[\log n-\lfloor nf(\varepsilon_1)\rfloor p\biggr] \nonumber \\
& \hspace{8pt} \geq 1-\exp\biggl[\log n-(nf(\varepsilon_1)-1) p\biggr] \nonumber \\
& \hspace{8pt} = 1-\exp\left[ -f(\varepsilon_1)np\left( 1-\frac{1}{f(\varepsilon_1)}\left[ \frac{\log n}{np}+\frac{1}{n} \right] \right) \right] \nonumber \\
& \hspace{8pt} = 1-\exp\biggl[-f(\varepsilon_1)np\bigl(1-o(1)\bigr)\biggr] & \text{(since } p=\omega\left(\frac{\log n}{n}\right) \text{)} \nonumber
% & \hspace{10pt} = 1-P_{\varepsilon_1}
\end{align}
which proves \eqref{eq:avshaperesult2cor}.
\end{proof}
Note that  if $\varepsilon_1$ is chosen very small in Corollary~\ref{cor:averageshape}, the inequalities will be empty statements in that the bounds on the probabilities will be negative. For example, since $\frac{\log n}{np}\rightarrow 0$ as $n\rightarrow\infty$, we have $\frac{1}{np}<\frac{\log n}{np}<\frac{1}{2}$ for $n>N$ for some $N>e$. It is easy to verify\footnote
{
Let us denote $x:=f(\varepsilon_1)=\frac{\varepsilon_1}{2+\varepsilon_1}$ and suppose $x>\frac{1+\log n}{np}$.
For the first probability limit we have (for $n>N$)
$$
2nxe^{-npx} < 2n\frac{1+\log n}{np}\frac{1}{ne}=\frac{2}{e}\frac{\log n}{np}<\frac{2}{e}\cdot\frac{1}{2}<1,
$$
since the function $x\mapsto 2nxe^{-npx}$ is decreasing for $x>\frac{1}{np}$ (hence also for $x>\frac{1+\log n}{np}$). For the second probability limit we have (for all $n$)
$$
\exp\bigl[\log n-(nf(\varepsilon_1)-1) p\bigr] = \frac{n}{e^{(nx-1)p}}=\frac{ne^p}{e^{npx}}<\frac{ne}{e^{npx}}<\frac{ne}{e^{1+\log n}}=1.
$$
}
that for these $n>N$, the probability bounds
\[
1-2nf(\varepsilon_1)\exp\left(-npf(\varepsilon_1)\right) \quad \text{and} \quad 1-\exp\bigl[\log n-(nf(\varepsilon_1)-1) p\bigr]
\]
in Corollary~\ref{cor:averageshape} are nonnegative if
\[
f(\varepsilon_1) > \frac{1+\log n}{np}.
\]
in which case $\varepsilon_1>\sqrt{\frac{\log n}{np}}$.
\section{Proof of the limit shape result}
Below follows the proof of Theorem~\ref{thm:main}.
\begin{proof}
Let $0<\varepsilon_1<\varepsilon$. Define the function $f(x)=\frac{x^2}{2+x}$, and  let us introduce the notation $P_{\varepsilon_1}$ as a shorthand for ${\exp\bigl[-f(\varepsilon_1)np(1-o(1))\bigl]}$. Also, as in the proof of Corollary \ref{cor:averageshape}, let us abbreviate $p=p_n$. 

By Corollary \ref{cor:averageshape}, $P_{\varepsilon_1}$ bounds both the probability that any pile in the first regime deviates from its expected size with more than $\varepsilon_1 np$ and the probability that any bowl in the second regime is nonempty. Now we will investigate how these pile sizes relate to the limiting function. (Note that we still have not rescaled the pile sizes.) Again, we treat two regimes of $k$ separately.

\underline{Regime 1}: $k\leq m$.
By the uniform convergence \eqref{eq:uniformconv2} in Observation~\ref{prop:uniformconv}, for any $\varepsilon_2>0$ and for sufficiently small $p$, we have
\begin{equation}
\label{eq:usecor}
|np(1-p)^{k-1}-npe^{-p(k-1)}|<np\varepsilon_2
\end{equation}
for all $k\geq 1$.
Since we chose $\varepsilon_1<\varepsilon$, we can choose $\varepsilon_2>0$ such that $\varepsilon_1+\varepsilon_2<\varepsilon$. Then, for sufficiently small $p$ (i.e., for sufficiently large $n$) we have with probability at least $1-P_{\varepsilon_1}$ (initially using the triangle inequality),
\begin{align}
\label{eq:regime1final}
& |\alpha_k-npe^{-p(k-1)}| \leq |\alpha_k-E\alpha_k|+|E\alpha_k-npe^{-p(k-1)}| \nonumber \\
&\hspace{20pt} = |\alpha_k-E\alpha_k|+|np(1-p)^{k-1}-npe^{-p(k-1)}| & \text{(by }\eqref{eq:Ebetak}\text{)} \nonumber \\
&\hspace{20pt} \leq np\varepsilon_1 + np\varepsilon_2 & \text{(by }\eqref{eq:regime1end},\eqref{eq:usecor}\text{)} \nonumber \\
&\hspace{20pt} = np(\varepsilon_1+\varepsilon_2).
\end{align}

\underline{Regime 2}: $k>m=\lfloor nf(\varepsilon_1) \rfloor$. In this regime, for sufficiently large $n$ (i.e., for sufficiently small $p$) with probability at least $1-P_{\varepsilon_1}$, we have
\begin{align}
\label{eq:regime2final}
|\alpha_k-npe^{-p(k-1)}| & \leq |\alpha_k|+|npe^{-p(k-1)}| \nonumber \\
& = 0 + npe^{-p(k-1)} & \text{(by }\eqref{eq:avshaperesult2cor}\text{)} \nonumber \\
& < npe^{-p(nf(\varepsilon_1)-1)} & \text{(since }k\in\mathbb{Z},k>nf(\varepsilon_1)\text{)} \nonumber \\
& = npe^p e^{-pnf(\varepsilon_1)} \nonumber \\
& \leq np(\varepsilon_1+\varepsilon_2),
\end{align}
where the last inequality follows from the assumption $p=\omega\left(\frac{\log n}{n}\right)$. (Namely, $np\rightarrow \infty$ as $n\rightarrow\infty$, so $e^p e^{-pnf(\varepsilon_1)}$ is eventually smaller than $\varepsilon_1+\varepsilon_2$.)

Combining \eqref{eq:regime1final} and \eqref{eq:regime2final}  and dividing by $np$ we obtain
\begin{equation}
\label{eq:endk}
P\left(\forall k>0: \left|\frac{\alpha_k}{np}-e^{-p(k-1)}\right|\leq \varepsilon_1+\varepsilon_2 \right) \geq 
1-P_{\varepsilon_1}.
\end{equation}

Now, we would like to say something about the $(1/p)$-rescaled boundary function $\rescaled{1/p}{\alpha}(x)=\frac{1}{np}\alpha_{\lfloor x/p \rfloor+1}$ for all real $x\geq 0$, and not only about the pile sizes $\alpha_k$ themselves for integers $k>0$. To this end, first define $\varepsilon_3>0$ such that $\varepsilon=\varepsilon_1+\varepsilon_2+\varepsilon_3$. Then, using the triangle inequality, we have for sufficiently large $n$ (i.e., for sufficiently small $p$) with probability at least $1-P_{\varepsilon_1}$,
\begin{align*}
\left|\frac{1}{np}\alpha_{\lfloor x/p\rfloor +1}-e^{-x}\right|
\leq \left|\frac{1}{np}\alpha_{\lfloor x/p\rfloor +1}-e^{-\lfloor x/p \rfloor p}\right|+\left|e^{-\lfloor x/p\rfloor p}-e^{-x}\right|
\end{align*}
\vspace{-10pt}
\begin{align}
& \hspace{30pt} \leq (\varepsilon_1+\varepsilon_2) + (e^{-\lfloor x/p\rfloor p}-e^{-x}) && \text{(by \eqref{eq:endk} and since }\lfloor x/p \rfloor p\leq x\text{)} \nonumber \\
& \hspace{30pt} \leq \varepsilon_1+\varepsilon_2+e^{-(x-p)}-e^{-x} && \text{(since }\lfloor x/p \rfloor p> x-p \text{)} \nonumber \\
& \hspace{30pt} \leq \varepsilon_1 + \varepsilon_2 + (e^p-1) && \text{(since }e^{-x}\leq 1 \text{)} \nonumber \\
& \hspace{30pt} \leq \varepsilon_1+\varepsilon_2+\varepsilon_3 && \text{(since }p\rightarrow 0\text{ as }n\rightarrow\infty \text{)} \nonumber \\
& \hspace{30pt} = \varepsilon \label{eq:thm1slut}.
\end{align}
Since the bound $\varepsilon$ in \eqref{eq:thm1slut} is independent of $x$, we conclude that it holds for all $x\geq 0$, and since also $\frac{1}{np}\alpha_{\lfloor x/p \rfloor+1}=\rescaled{1/p}{\alpha}(x)$, we can write \eqref{eq:thm1slut} as
$$
\Phi:=P\left( \|\rescaled{1/p}{\alpha} - g\|_{\infty} \leq \varepsilon \right) \geq
1-P_{\varepsilon_1}.
$$
In other words, $\liminf\limits_{n\rightarrow\infty}\frac{-\log(1-\Phi)}{np} \geq f(\varepsilon_1)$ and, since this holds for any $\varepsilon_1<\varepsilon$, the continuity of the function $f$ yields that $\liminf\limits_{n\rightarrow\infty} \frac{-\log(1-
\Phi)}{np} \geq f(\varepsilon)$. Thus,
$$
P\left( \|\rescaled{1/p}{\alpha} - g\|_{\infty} \leq \varepsilon \right) \geq 1-P_{\varepsilon}=
1-\exp\biggl[-f(\varepsilon)np\bigl(1-o(1)\bigr)\biggr].
$$
An application of Lemma~\ref{lem:reordering} concludes the proof of Theorem~\ref{thm:main}.
\end{proof}

\section{Discussion}
Popov initiated the study of stochastic versions of Bulgarian solitaire \cite{Popov}. Here we compare our work with his. Like the original deterministic game, Popov's stochastic version was ``pile-based'' in that each old pile was independently picked with a fixed probability $p$ to release a card to the new pile.  In contrast, we here studied a ``card-based'' version, where every single card was independently picked with the same probability $p$. This modification radically changed the outcome of the process from a triangular limit shape to an exponential limit shape. From an analytical point of view, a crucial difference is that the card-based version allowed us to study each pile independently of previously formed piles. To capitalise on this feature we represented game configurations by weak integer compositions. We proved asymptotic results about the shape of these compositions and thanks to Lemma~\ref{lem:reordering} these results would hold also for integer partitions obtained by sorting the parts of compositions.

By formulating our results in terms of a limit shape of integer partitions we connected to an important literature by Vershik and others. However, in that literature the starting point is typically a given probability distribution on integer partitions, see \cite{VershikStatMech}. Here the starting point was instead a stochastic process on Young diagrams. Often a stochastic process defines a unique stationary probability distribution, which in turn will define a limit shape. This was the case in our Corollary~\ref{thm:statmeasure} as well as in previous work on birth-and-death processes on Young diagrams by Eriksson and Sj\"ostrand \cite{Eriksson2012575}. Although Popov did not connect his result to Vershik's theory of limit shapes of integer partitions, his result too can be interpreted as a triangular limit shape of this kind being obtained under the stationary probability distribution of the pile-based stochastic Bulgarian solitaire. 

A novel feature of the present work is that we proved a limit shape result that is genuinely about the process instead of just using its stationary distribution: Theorem~\ref{thm:main} gives a relation between how many moves are played and the probability that this sequence of moves will end up close to an exponential shape. Theorem~\ref{thm:averageshape} specifically says that after exactly $m$ steps in the process, the number of cards in the piles are close to their expected value with a (high) probability that depends on $m$. This is kind of question could be asked for any stochastic process on the partitions of a fixed integer $n$, such as those in \cite{Eriksson2012575} and \cite{Popov}.

\bibliographystyle{amsplain}
\bibliography{limitshape}	% expects file "limitshape.bib"
\end{document}